\documentclass[a4paper,12pt]{article}
\usepackage{dcolumn}% Align table columns on decimal point
\usepackage{bm}% bold math
\usepackage[dvips]{graphicx}
\usepackage{amsmath}
\usepackage{hyperref}
\usepackage{amssymb}%Nik
\usepackage{amsfonts}
\usepackage{ucs}
\usepackage[utf8x]{inputenc}
\usepackage[T1]{fontenc}
\usepackage{url}
\usepackage{psfrag}
\usepackage{dsfont}
\usepackage{color}

\newtheorem{theorem}{\textbf{Theorem}}

\begin{document}

\title{Approximation of a compressible Navier-Stokes system by  non-linear acoustical models}

\author{ANNA ROZANOVA-PIERRAT\footnote{Laboratoire Math\'ematiques et Informatique Pour la Complexit\'e et les Syst\`emes,
Centrale Sup\'elec, Universit\'e Paris-Saclay, Grande Voie des Vignes,
Ch\^atenay-Malabry, France
anna.rozanova-pierrat@centralesupelec.fr}}
%\subtitle{{\small Document prepared by Anna Rozanova-Pierrat\footnotemark}}
% \address{Laboratory Applied Mathematics and Systems, Ecole Centrale Paris, Grande Voie des Vignes\\
% Ch\^atenay-Malabry, France\\
% anna.rozanova-pierrat@ecp.fr}

\maketitle

%\index{Rozanova-Pierrat, A.V.}                              % write this for each author
% \index{Author2, I.I.}                              % to generate the index
% \index{Coauthor, I.I.}                             %

\begin{abstract}
     We analyse the existing derivation of the models of non-linear acoustics such as the Kuznetsov equation, the NPE equation and the KZK equation. The technique of introducing a corrector in the derivation \emph{ansatz} allows to consider the solutions of these equations as approximations of the solution of the initial system (a compressible Navier-Stokes/Euler system). The validation of the approximation \emph{ansatz} is given for the KZK equation case.
   \end{abstract}

\section{Introduction}
There is a renewed interest in the study of wave propagation, in particular
because of recent applications to ultrasound imaging (i.e. HIFU) or technical and
medical applications such as lithotripsy or thermotherapy. Such new techniques
rely heavily on the ability to model accurately  the nonlinear propagation of a
 finite-amplitude sound pulse in thermo-viscous elastic media.
%Body of the paper divided into sections and subsections.

We analyse the derivation of different models of  non-linear acoustics such as the Kuznetzov~\cite{Kuznetsov}, the Nonlinear Progressive wave Equation (NPE)~\cite{MCDKuper} and the Khokhlov-Zabolotskaya-Kuznetzov  (KZK)~\cite{KhZabKuzn} equations which are perturbative and paraxial approximations of small perturbations around a given state of  a compressible nonlinear isentropic Navier-Stokes (for  viscous media) and Euler (for the non-viscous case) systems.
The direct derivation shows that the Kuznetzov equation is the first order approximation of the Navier-Stokes system, the KZK and NPE equations are the first order approximations of the Kuznetzov equation and the second order approximations of the Navier-Stokes system.
In addition,  the NPE equation can be considered as an approximation of the KZK equation.

%Therefore, using the  well-posedness results for the Navier-Stokes/Euler systems and for the KZK equation~\cite{ARPAnal}, we find the domains, where we validate the approximation of the exact solution of the Navier-Stokes/Euler systems by the solution of the KZK equation.

%KZK can be also obtained by other types of derivation~\cite{DJMR,Sanchez,Tex}.
To be able to validate the approximation of the exact solution of the Navier-Stokes/Euler systems by the solution of the Kuznetsov/KZK/NPE equation,  we need to ensure that the derivation of our model, the Kuznetsov/KZK/NPE equation, allows us to  reconstruct
the solution of the initial Navier-Stokes system from the solution of the Kuznetsov/KZK/NPE equation.
In this aim, following the ideas of Refs.~\cite{ARPDer,Anna}, we modify the initial physical derivation, given in Refs.~\cite{Kuznetsov,KhZabKuzn} for the KZK and the Kuznetsov equations and given in Ref.~\cite{MCDKuper} for the NPE equation, introducing a corrector function in the derivation \textit{ansatz}.

We also improve
 the validation of the KZK-approximation for the non-viscous and viscous
cases obtained in Ref.~\cite{ARPDer}, by the precision of the speed order of divergence between the solutions of the approximate and the exact systems.

 Let us introduce some notations used throughout the paper.
%The operator $\nabla.$ denotes the divergence operator.
For   a positive fixed  small enough real number $\epsilon$,  we suppose
that $\mathbb{R}_+$ consists of classes, which are characterized by the
power of $\epsilon$:
$$\ldots,\epsilon^2,\ldots,\epsilon,\ldots,\sqrt{\epsilon},
\ldots, \epsilon^0 =1,\ldots,  \frac 1{\epsilon},\ldots, \frac 1{\epsilon^2},\ldots
   $$   $O(1)$  denotes the class of constants.
 %If $U$ is a matrix, $U^T$ is the transpose of $U$.
 
\section{Approximation of the hydro-dynamic system by  an isentropic Navier-Stokes system}
  We start from the Navier-Stokes system in $\mathbb{R}^n$:
\begin{gather} 
\partial_t \rho  +\operatorname{div}( \rho  \mathbf{u})=0,\label{nss1}\\
\rho[\partial_t
\mathbf{u}+(\mathbf{u} \cdot\nabla) \,
\mathbf{u}] = -\nabla p+ \beta \nabla\operatorname{div} \mathbf{u},\label{nss2}\\
\rho T \left[\partial_t S+(\mathbf{u} \cdot \nabla) S \right]=\kappa \triangle T+\zeta (\operatorname{div}\mathbf{u})^2\nonumber\\
+\frac{\eta}{2}\left(\partial_{x_k}u_i+\partial_{x_i} u_k-\frac{2}{3}\delta_{ik}\partial_{x_i}u_i \right)^2,\label{nss3}\\
p=p(\rho, S) ,\label{nss4}%\\
%  \rho[\partial_t
% u+(u\cdot\nabla) u] = -\nabla p(\rho, S) +b \Delta u\,,\label{nss2}
\end{gather}
where $S$ is
the entropy and the state law $p=p(\rho, S)$ is the pressure.  The density $\rho$, the velocity $\mathbf{u}$, the temperature $T$ and the entropy are unknown functions in the system~(\ref{nss1})--(\ref{nss4}). The coefficients $\beta$, $\kappa$ and $\eta$ are constant viscosity coefficients.  

First, we assume that the temperature $T$ and the entropy $S$ have
 small increments $T=T_0+\epsilon\tilde{T}$ and $S=S_0+\epsilon^2\tilde{S}$. With the
hypothesis of potential motion, we introduce  constant states
$$\rho=\rho_0,\quad \mathbf{u}=\mathbf{u_0}.$$
  Next, we assume that the  density fluctuations (around the
constant state $\rho_0$) and the velocity fluctuations (around $\mathbf{u_0}$, which can be
taken equal to zero using a  Galilean transformation), are of the
same order of
  $\epsilon$:
  \begin{equation}\label{ap1}
  \rho_\epsilon=\rho_0+ \epsilon  \tilde{\rho}_{\epsilon}\,,\quad \mathbf{u_\epsilon}=\epsilon\mathbf{\tilde
u_\epsilon}\,,
  \end{equation}
where $\epsilon $ is a dimensionless parameter which characterizes the
smallness of the perturbation. For instance, in water with an initial
power of the order of $0.3 \, \mathrm{W} /\mathrm{cm}^2$  $\epsilon$ is equal to $10^{-5}$. We also suppose that all viscosity coefficients, for instance, $\beta,$ $\zeta$, $\eta$ and $\kappa$, are small of the order $\epsilon$:
$$ \beta=\epsilon \tilde{\beta}.$$
     Using the transport heat
equation up to the terms of the order of $\epsilon^3$
$$\epsilon^2\rho_0 T_0 \partial_t
\tilde{S}=\epsilon^2\tilde{\kappa} \triangle \tilde{T}+O(\epsilon^3),$$ %where $\kappa$  is a dissipation coefficient,
the approximate state equation 
$$
 p=p_0+c^2 \epsilon\tilde{\rho}_{\epsilon} + \frac{1}{2}\left( \partial^2_\rho p\right)_S
{\epsilon}^2{\tilde{\rho}_{\epsilon}}^2 +\left(\partial_S p
\right)_{\rho} {\epsilon}^2\tilde{S}+O(\epsilon^3)
$$
(where the notation $\left(\cdot\right)_S$ means that the expression
in brackets is constant in $S$), can be replaced~\cite{KhZabKuzn,Makarov,hamBlack}
%,
%due to the following relation $$S-S_0=-\frac{\kappa}{T_0}\left(\partial_p
%T\right)_S \operatorname{div} u_\epsilon ,$$ 
by
\begin{equation}
p=p_0+c^2 \epsilon{\tilde{\rho}}_{\epsilon} +\frac{(\gamma-1) c^2}{2 \rho_0}{\epsilon}^2
{{\tilde{\rho}}_{\epsilon}}^2
- \epsilon\tilde{\kappa}\left( \frac{1}{C_v} -\frac{1}{C_p}
\right) \nabla. u_{\epsilon}+O(\epsilon^3).\label{s2}
\end{equation}
Here $\gamma=C_p/C_v$ denotes the ratio of the heat capacities at  constant pressure and at  constant volume respectively.
System~(\ref{nss1})--(\ref{nss4}) becomes an isentropic system
\begin{gather} 
\partial_t \rho_\epsilon  +\operatorname{div}( \rho_\epsilon  \mathbf{u_\epsilon})=0\,,\label{NSi1}\\
 \rho_\epsilon[\partial_t
\mathbf{u_\epsilon}+(\mathbf{u_\epsilon} \cdot\nabla) \,
\mathbf{u_\epsilon}] = -\nabla p(\rho_\epsilon)
+\epsilon \nu \Delta \mathbf{u_\epsilon}\,,\label{NSi2}
\end{gather}
with the approximate state equation 
\begin{equation}\label{press}
p(\rho_{\epsilon})=p_0+c^2 (\rho_{\epsilon}-\rho_0) +\frac{(\gamma-1) c^2}{2
\rho_0}(\rho_{\epsilon}-\rho_0)^2\quad (O(\epsilon^3))
\end{equation}
and with a small enough and positive  viscosity coefficient:
$$\epsilon \nu =\beta+\kappa\left( \frac{1}{C_v} -\frac{1}{C_p} \right).$$
\section{Perturbative approach: Kuznetsov equation}
% \subsection{Linearised Navier-Stokes system: the wave equation}
% Let us start with a simple linearization of the isentropic Navier-Stokes system~(\ref{NSi1})--(\ref{NSi2}) with linear approximation of the state equation~(\ref{press}) $p(\rho_{\epsilon})=p_0+c^2 (\rho_{\epsilon}-\rho_0)+O(\epsilon^2)$, taking in consideration relations~(\ref{ap1}):
% \begin{align}
%      &\partial_t \rho_\epsilon  +\operatorname{div}( \rho_\epsilon  u_\epsilon)=\epsilon\left[\partial_t \rho_1+\rho_0 \operatorname{div}\tilde{u}\right]+\epsilon^2 \operatorname{div}(\rho_1\tilde{u}), \label{ns1Lin}\\
% &\rho_\epsilon[\partial_t
% u_\epsilon+(u_\epsilon \cdot\nabla) \,
% u_\epsilon] +\nabla p(\rho_\epsilon)
% -\epsilon \nu \Delta u_\epsilon=\nonumber\\
% &\epsilon\left[\rho_0\partial_t \tilde{u}+c^2\nabla \rho_1\right]+\epsilon^2\left[\rho_1 \partial_t \tilde{u}+\rho_0(\tilde{u}\cdot \nabla)\tilde{u}-\nu \Delta \tilde{u}\right]\nonumber\\
% &+\epsilon^3 \rho_1
% (\tilde{u}\cdot \nabla)\tilde{u}.       \label{ns2Lin}   
%               \end{align}
% 
% 
% 
% \subsection{Kuznetsov equation}
First derived by Kuznetsov~\cite{Kuznetsov} from the isentropic Navier-Stokes system~(\ref{NSi1})--(\ref{press}), the Kuznetsov equation
\begin{equation}\label{KuzEq}
 \partial^2_t \tilde{\phi} -c^2 \triangle \tilde{\phi}=\partial_t\left( (\nabla \tilde{\phi})^2+\frac{\gamma-1}{2c^2}(\partial_t \tilde{\phi})^2+\frac{\epsilon\nu}{\rho_0}\Delta \tilde{\phi}\right),
\end{equation}
written for the velocity potential $$\mathbf{u}(\mathbf{x},t)=-\nabla \tilde{\phi}(\mathbf{x},t),\quad \mathbf{x}\in \mathbb{R}^n, \quad t\in \mathbb{R}^+,$$ was latter derived by other methods and was discussed by a lot of authors (see for examples~\cite{hamBlack,Jordan}).

Here we focus on the introduction of the corrector $\epsilon^2\rho_2$ in the \emph{ansatz} of Kuznetsov
\begin{gather}
\rho_\epsilon(\mathbf{x}, t)=\rho_0+\epsilon\rho_1(\mathbf{x},t)+\epsilon^2\rho_2(\mathbf{x},t)\label{rhoKuz}\\
\mathbf{u_\epsilon} (\mathbf{x}, t)=-\epsilon\nabla \phi(\mathbf{x},t),\label{uKuz}
\end{gather}
which allows to open the question about the approximation between the exact solution of the isentropic Navier-Stokes system~(\ref{NSi1})--(\ref{press}) and its approximation given by the solution of the Kuznetsov equation, as it was done for the KZK equation in~\cite{ARPDer}.

Puting expressions for the density and the velocity~(\ref{rhoKuz})--(\ref{uKuz}) into the  isentropic Navier-Stokes system~(\ref{NSi1})--(\ref{press}), we directely obtain
 \begin{align}
     &\partial_t \rho_\epsilon  +\operatorname{div}( \rho_\epsilon  \mathbf{u_\epsilon})=\epsilon\frac{\rho_0}{c^2}\left[\partial^2_t \phi -c^2 \triangle \phi-\right.\nonumber\\
     &\left.\epsilon\partial_t\left( (\nabla \phi)^2+\frac{\gamma-1}{2c^2}(\partial_t \phi)^2+\frac{\nu}{\rho_0}\Delta \phi\right)\right]+O(\epsilon^3), \label{ns1K}\\
&\rho_\epsilon[\partial_t
\mathbf{u_\epsilon}+(\mathbf{u_\epsilon} \cdot\nabla) \,
\mathbf{u_\epsilon}] +\nabla p(\rho_\epsilon)
-\epsilon \nu \Delta \mathbf{u_\epsilon}=\nonumber\\
&\epsilon\nabla\left[\rho_1 -\frac{\rho_0}{c^2}\partial_t \phi\right]+\epsilon^2\nabla\left[c^2 \rho_2+\frac{\rho_0(\gamma-2)}{2c^2} (\partial_t \phi)^2\right.\nonumber\\
&\left.+\frac{\rho_0}{2}(\nabla \phi)^2+\nu\Delta \phi\right]+O(\epsilon^3).       \label{ns2K}   
              \end{align}
We see that the Kuznetsov equation
\begin{equation}\label{KuzEqA}
 \partial^2_t \phi -c^2 \triangle \phi=\epsilon\partial_t\left( (\nabla \phi)^2+\frac{\gamma-1}{2c^2}(\partial_t \phi)^2+\frac{\nu}{\rho_0}\Delta \phi\right),
\end{equation}
is the first order approximation, obtained from the equation of mass conservation up to the terms $O(\epsilon^3)$ with the relations for the density perturbations, found from the momentum conservation also up to the terms $O(\epsilon^3)$ with the help of the Sommerfeld radiation boundary condition at infinity:
\begin{gather}%\begin{equation}
 \rho_1(\mathbf{x},t) =\frac{\rho_0}{c^2}\partial_t \phi(\mathbf{x},t),\label{rho1K}\\ \rho_2(\mathbf{x},t)=-\frac{\rho_0(\gamma+2)}{2c^4} (\partial_t \phi)^2-\frac{\rho_0}{2c^2 }(\nabla \phi)^2-\frac{\nu}{c^2}\Delta \phi.\label{rho2K}
\end{gather}
Since initially, we consider the state equation for the pressure $p$ up to the terms of the order of $\epsilon^3$, we  conclude that the \emph{ansatz} of the Kuznetsov equation gives the optimal approximation error of the same order.

Let us also notice, as it was  originaraily mentionned  by Kuznetsov, that the Kuznetsov equation~(\ref{KuzEqA}) contains terms of different orders, and hence, it is a wave equation with small size non-linear perturbations  $\partial_t(\nabla \phi)^2$, $\partial_t(\partial_t \phi)^2$ and  viscosity term $\partial_t\Delta \phi$.
A way to obtain an approximate equation containning all terms of the same  order without modification  of the order of remainder terms is to perform a paraxial approximation, which we introduce in the next section. This time   the approximation becomes the second order approximation and will be given by the KZK equation. 
\section{Paraxial approximation}
\subsection{KZK equation}
In the present Section we focus on the derivation of 
 the KZK equation~(\ref{KZKI})
     in non-linear media using the following acoustical properties of beam's propagation
\begin{enumerate}
  \item The beams are concentrated near the $x_1$-axis ;
  \item The beams propagate along the $x_1$-direction;
  \item The beams are generated either by an initial condition or by a forcing term on the boundary $x_1=0$.
  \end{enumerate}
 It is assumed that the variation of  beam's
propagation in the direction
$$
\mathbf{x'}=(x_2,x_3,\ldots, x_n)
$$
perpendicular to the $x_1$-axis is much larger than its variation
along the $x_1$-axis, \textit{i.e.} we suppose that the beam has the form $U(t-x_1/c,\epsilon x_1,\sqrt{\epsilon}\mathbf{x'})$. The first argument $t-x_1/c$ describes the wave propagation in time along the $x_1$-axis with the sound speed $c$, two last arguments $\epsilon x_1$ and $\sqrt{\epsilon}\mathbf{x'}$ describe respectively   the speed of the deformation of the wave along the $x_1$-axis and along the $\mathbf{x'}$-axis. We remark that $\epsilon\ll 1$ and consequently, $\epsilon\ll \sqrt{\epsilon}$.

\begin{figure}[t!]
\begin{center}
                   \psfrag{a}{$x_1$}\psfrag{b}{$\mathbf{x'}$}\psfrag{c}{$t$}\psfrag{NS}{\small{Navier-Stokes/}}\psfrag{E}{\small{ Euler $(x_1,\mathbf{x'},t)$}}\psfrag{a1}{$z=\epsilon x_1$}\psfrag{b1}{$y=\sqrt{\epsilon} \mathbf{x'}$}\psfrag{c1}{$\tau=t-\frac{x_1}{c}$}\psfrag{KZK}{\small{KZK $(\tau,z,\mathbf{y})$}}
  \includegraphics[width=.43\textwidth]{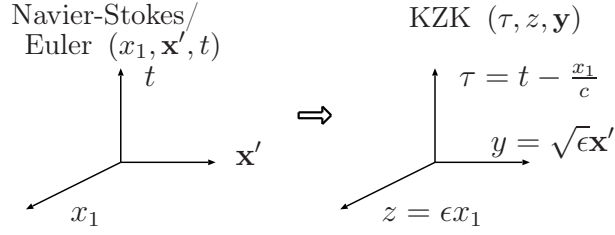}
                  \end{center}
   \caption{Paraxial change of variables for the profiles $U(t-x_1/c,\epsilon x_1,\sqrt{\epsilon}\mathbf{x'})$.}\label{fig1}
% \end{figure}
%\includegraphics[width=0.9\linewidth]{./Base.eps}
\end{figure}

We notice that  if we perform the paraxial change of variables (see Fig.~\ref{fig1}),
the wave operator $\partial_t^2-c^2\Delta$ becomes
$$\partial_t^2-c^2\Delta= \epsilon\left[2c\partial_{\tau z}^2-c^2\Delta_\mathbf{y}\right]-\epsilon^2c^2\partial_z^2 .$$
Therefore, if we suppose that the velocity potential $\phi(\mathbf{x},t)=\Phi(t-x_1/c,\epsilon x_1,\sqrt{\epsilon}\mathbf{x'})$, we directely obtain from the Kuznetsov equation~(\ref{KuzEqA}) (see also~\cite{Kuznetsov}) that
\begin{align}
 &\partial^2_t \phi -c^2 \triangle \phi-\epsilon\partial_t\left( (\nabla \phi)^2+\frac{\gamma-1}{2c^2}(\partial_t \phi)^2+\frac{\nu}{\rho_0}\Delta \phi\right)\nonumber\\
&=\epsilon\left[2c\partial_{\tau z}^2\Phi-\frac{\gamma+1}{4c^2}c^2\partial_\tau(\partial_\tau \Phi)^2\right.\nonumber\\
&\left.-\frac{\nu}{\rho_0c^2}\partial^3_\tau\Phi-\Delta_\mathbf{y} \Phi\right]+O(\epsilon^2).\label{KZKpoten}
\end{align}
%Now using relation~(\ref{rho1K}), we find that

Therefore, returning to the derivation of the Kuznetsov equation, after the paraxial approximation of  $\phi$, $\rho_1$ and $\rho_2$ with profiles $\Phi$, $I$ and $J$ 
\begin{align*}
&\mathbf{u_\epsilon}(\mathbf{x},t)=-\epsilon\left(-\frac{1}{c}\partial_\tau\Phi+\epsilon \partial_z \Phi; \sqrt{\epsilon}\nabla_\mathbf{y} \Phi\right)(\tau,z,\mathbf{y})\\
&\rho_1(\mathbf{x},t)=I(\tau,z,\mathbf{y})=\frac{\rho_0}{c^2}\partial_\tau \Phi(\tau,z,\mathbf{y}), \\
 &\rho_2(\mathbf{x},t)=J(\tau,z,\mathbf{y})=\\&-\frac{(\gamma-1)\rho_0}{2c^4}(\partial_\tau \Phi)^2-\frac{\nu }{c^4} \partial^2_\tau\Phi+O(\epsilon),
 \end{align*}
 we find that the right-hand $\epsilon$-order terms in Eq.~(\ref{KZKpoten}) is exactly the KZK equation, originally written in Ref.~\cite{KhZabKuzn} for the (first) perturbation $I$ of the density $\rho_\epsilon$:
 \begin{equation}\label{KZKI}
  c\partial^2_{\tau z} I -\frac{(\gamma+1)}{4\rho_0}\partial_\tau^2
I^2-\frac{\nu}{2 c^2\rho_0}\partial^3_\tau I-\frac{c^2}2 \Delta_\mathbf{y}
I=0.
 \end{equation}
We notice that this model still contains terms describing the wave propagation $\partial^2_{\tau z} I$,  the non-linearity  $ \partial_\tau^2
I^2$ and the viscosity effects $\partial^3_\tau I$ of the medium, as the Kuznetsov equation and adds a diffraction effects by the tranversal laplacian $ \Delta_\mathbf{y} I$.

 In addition, performing the paraxial approximation in the right-hand side of equations~(\ref{ns1K})--(\ref{ns2K}), we obtain that the KZK equation is the second order approximation of the isentropic Navier-Stokes system up to term of $O(\epsilon^3)$. %:% respectively becomes
%  \begin{align*}
%      &\partial_t \rho_\epsilon  +\operatorname{div}( \rho_\epsilon  u_\epsilon)=\epsilon^2\frac{\rho_0}{c^2}\left[2c\partial_{\tau z}^2\Phi-\frac{\gamma+1}{4c^2}c^2\partial_\tau(\partial_\tau \Phi)^2\right.\\
% &\left.-\frac{\nu}{\rho_0c^2}^2\partial^3_\tau\Phi-\Delta_y \Phi\right]+O(\epsilon^3),\\% \label{ns1K}\\
% &\rho_\epsilon[\partial_t
% u_\epsilon+(u_\epsilon \cdot\nabla) \,
% u_\epsilon] +\nabla p(\rho_\epsilon)
% -\epsilon \nu \Delta u_\epsilon=\nonumber\\
% &\epsilon\nabla\left[\rho_1 -\frac{\rho_0}{c^2}\partial_t \phi\right]+\epsilon^2\nabla\left[c^2 \rho_2+\right.\nonumber\\
% &\left.\frac{\rho_0(\gamma+2)}{2c^2} (\partial_t \phi)^2+\frac{\rho_0}{2}(\nabla \phi)^2+\nu\Delta \phi\right]+O(\epsilon^3).      % \label{ns2K}   
%               \end{align*}
In this sense, since the entropy and the pressure are approximated up to terms of the order of $\epsilon^3$, the Kuznetsov-type \emph{ansatz} (for the Kuznetsov or the KZK equations) is optimal, as the equations of the Navier-Stokes system also approximated up to $O(\epsilon^3)$-terms.
For instance, the \emph{ansatz} initially proposed by Khokhlov and Zabolotskaya~\cite{KhZabKuzn} to derive the KZK equation, corrected with $\epsilon^2 v_1$~\cite{ARPDer} for the velocity perturbation along the propagation axis,
 \begin{align*}
   &  \rho_\epsilon(x_1,\mathbf{x'},t)= \rho_0+\epsilon I(t-\frac{x_1}{c},\epsilon x_1, \sqrt{\epsilon} \mathbf{x'})\,,\\%\label{scal1}\\
   &\mathbf{u_\epsilon}(x_1,\mathbf{x'},t)=\epsilon(v+\epsilon v_1 ;\sqrt{\epsilon}
   \mathbf{w})(t-\frac{x_1}{c},\epsilon x_1, \sqrt{\epsilon} \mathbf{x'})%\label{scal2}\,.
  \end{align*}
is not optimal since the equation of momentum in tranversal direction keeps the non-zero terms of the order of $\epsilon^\frac{5}{2}$~\cite{ARPDer}.

\subsection{NPE equation}
\begin{figure}[t!]
\begin{center}
                   \psfrag{a}{$x_1$}\psfrag{b}{$\mathbf{x'}$}\psfrag{c}{$t$}\psfrag{NS}{\small{Navier-Stokes/}}\psfrag{E}{\small{ Euler $(x_1,\mathbf{x'},t)$}}\psfrag{a1}{$z=x_1-ct$}\psfrag{b1}{$\mathbf{y}=\sqrt{\epsilon} \mathbf{x'}$}\psfrag{c1}{$\tau=\epsilon t$}\psfrag{NPE}{\small{NPE $(\tau,z,\mathbf{y})$}}
  \includegraphics[width=.43\textwidth]{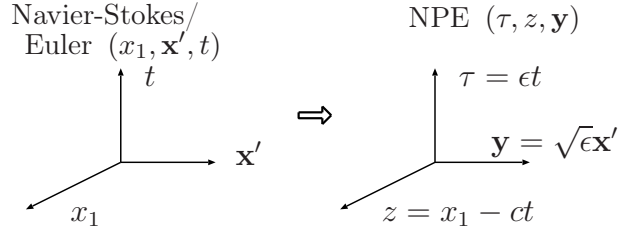}
                  \end{center}
   \caption{Paraxial change of variables for the profiles $U(\epsilon t,x_1-ct,\sqrt{\epsilon}\mathbf{x'})$.}\label{fig2}
% \end{figure}
%\includegraphics[width=0.9\linewidth]{./Base.eps}
\end{figure}

The NPE equation (Nonlinear Progressive wave Equation), initially derived by McDonald and Kuperman~\cite{MCDKuper},
gives another example of a paraxial approximation in the aim to describe short-time pulses and a long-range propagation (see Fig.~\ref{fig2}), for instance, in an ocean waveguide, where the refraction phenomena are important.
To compare to the KZK-\emph{ansatz}, the role of propagation distance and time was reversed~\cite{MCDKuper}:
$$z_{\mathrm{NPE}}=-c\tau_{\mathrm{KZK}},\quad \tau_{\mathrm{NPE}}=\epsilon\tau_{\mathrm{KZK}}+\frac{z_{\mathrm{KZK}}}{c}.$$
Consequently, from the KZK equation we directely have the NPE equation with the error $O(\epsilon)$:
\begin{align*}
 &c\partial^2_{\tau z} I -\frac{(\gamma+1)}{4\rho_0}\partial_\tau^2
I^2-\frac{\nu}{2 c^2\rho_0}\partial^3_\tau I-\frac{c^2}2 \Delta_\mathbf{y}
I=\\
&-c \partial^2_{\tau_{\mathrm{NPE}} z_{\mathrm{NPE}}} I -\frac{c^2(\gamma+1)}{4\rho_0}\partial_{z_{\mathrm{NPE}}}^2
I^2+\frac{c\nu}{2 \rho_0}\partial^3_{z_{\mathrm{NPE}}} I\\
&-\frac{c^2}2 \Delta_{\mathbf{y}_{\mathrm{NPE}}} I+O(\epsilon).
\end{align*}
The fact that the NPE equation is an approximation of the KZK equation does not allow to keep, by the analogy to the derivation of the KZK, the Kuznetsov-\emph{ansatz} of perturbations~(\ref{rhoKuz})--(\ref{uKuz}) just by introducing the new paraxial profiles $\Psi$ for $\phi$, $P_1$ for $\rho_1$ and $P_2$ for $\rho_2$. Indeed, if we do this,
the Kuznetsov equation, appeared in the conservation of  mass, gives the NPE equation for the potential profile $\Psi$ [compare with Eq.~(\ref{KZKpoten})]
\begin{align}
 &\partial^2_t \phi -c^2 \triangle \phi-\epsilon\partial_t\left( (\nabla \phi)^2+\frac{\gamma-1}{2c^2}(\partial_t \phi)^2+\frac{\nu}{\rho_0}\Delta \phi\right)\nonumber\\
&=\epsilon\left[-2c\partial_{\tau z}^2\Psi+\frac{\gamma+1}{2}c\partial_z(\partial_z \Psi)^2\right.\nonumber\\
&\left.+\frac{\nu c}{\rho_0}^2\partial^3_z\Psi-c^2\Delta_\mathbf{y} \Psi\right]+O(\epsilon^2),\label{NPEpoten}
\end{align}
but in the conservation of momentum,
we obtain that the corrector $P_1$ has a term of the order of $\epsilon$:
\begin{align*}
%&u_\epsilon(x,t)=-\epsilon\left(\partial_z\Psi; \sqrt{\epsilon}\nabla_y \Psi\right)(\tau,z,\mathbf{y})\\
&\rho_1(\mathbf{x},t)=P_1(\tau,z,\mathbf{y})=-\frac{\rho_0}{c}\partial_z \Psi +\epsilon\frac{\rho_0}{c^2} \partial_\tau\Psi, 
 \end{align*}
 what will not allow to keep equal to zero just the terms of the same order without any arrangement between the first and the second order terms.
 Thus we need to suppose that
 \begin{align*}
  &\mathbf{u_\epsilon}(\mathbf{x},t)=-\epsilon\nabla \phi(\mathbf{x},t)=-\epsilon\left(\partial_z\Psi; \sqrt{\epsilon}\nabla_\mathbf{y} \Psi\right)(\tau,z,\mathbf{y}),\\
  &\rho_\epsilon(\mathbf{x},t)=\rho_0+\epsilon P_1(\tau,z,\mathbf{y})+\epsilon^2 P_2(\tau,z,\mathbf{y}),
   \end{align*}
where
\begin{align*}
 &P_1(\tau,z,\mathbf{y})=\frac{\rho_0}{c}\partial_z \Psi(\tau,z,\mathbf{y}),\\
 &P_2(\tau,z,\mathbf{y})=\frac{\rho_0}{c^4}\partial_\tau \Psi-\frac{\rho_0(\gamma+3)}{2c^2}(\partial_z \Psi)^2-\frac{\nu}{c^2}\partial_z^2\Psi,
\end{align*}
to obtain the NPE equation for the profile of the potential
\begin{align}
& \partial^2_{\tau z} \Psi-\frac{\gamma+1}{4} \partial_z(\partial_z \Psi)^2 - \frac{\nu}{2 \rho_0} \partial_z^3\Psi+\frac{c}{2}\Delta_\mathbf{y} \Psi=0\label{NPE}
\end{align}
as the second order approximation of the isentropic Navier-Stokes system up to the terms of the order of $O(\epsilon^3)$.

\section{Approximation results}
We precise the approximation results for the KZK equation, given in Ref.~\cite{ARPDer}, by the evaluation of the size of the difference between the exact and the approximate solutions. 
As it was explained in Ref.~\cite{ARPDer}, the isentropic Euler system for  $\mathbf{\widetilde{U}_\epsilon}=(\rho_\epsilon,\rho_\epsilon \mathbf{u_\epsilon})$ and $\mathbf{F(\widetilde{U}_\epsilon)}=(\rho_\epsilon \mathbf{u_\epsilon}, \rho_\epsilon \mathbf{u^2_\epsilon} +p(\rho_\epsilon))^T$ can be written as a system of conservation laws
\begin{equation}\label{CCL}
 \partial_t \mathbf{\widetilde{U}_{\epsilon}} +\nabla .\mathbf{F(\widetilde{U}_{\epsilon})}=0.
\end{equation}
The KZK-ansatz allows to find from the solution $I$ of the KZK equation~(\ref{KZKI}) the correctors $v$, $\mathbf{w}$,$v_1$ and to obtain
for \begin{equation}\label{linU}
\mathbf{\overline
U_\epsilon}=(\overline{\rho}_\epsilon,\overline{\rho}_\epsilon\mathbf{\overline{u}_\epsilon}),
\end{equation}
with \begin{align*}
      &\overline{\rho}_\epsilon =\rho_0+\epsilon
I(t-\frac{x_1}c,\epsilon x_1,
\sqrt{\epsilon} \mathbf{x'}),\\
&\mathbf{\overline{u}_\epsilon}=\epsilon(v+\epsilon v_1,\sqrt{\epsilon} \mathbf{w}))(t-\frac{x_1}c,\epsilon x_1,
\sqrt{\epsilon} \mathbf{x'}),
      \end{align*}
the approximate system
\begin{equation}\label{CL}
 \partial_t \mathbf{\overline{U}_{\epsilon}} +\nabla .\mathbf{F(\overline{U}_{\epsilon})}=\epsilon^\frac{5}{2}\mathbf{R}.
\end{equation}

More precisely, for the non-viscous case, we have the  following theorem:
\begin{theorem}
 \label{Euler}
Let $I_0(\tau,0,\mathbf{y})\in H^{s'}(\mathbb{R}^n)$, $s'> [\frac n2]+ 5$ be the initial data for the KZK equation~(\ref{KZKI}) $L$-periodic and with mean value zero with respect to $\tau$.
Then there exists a unique solution $I$ of the KZK equation
     such that
    \begin{itemize}
        \item $I(\tau,z,\mathbf{y})$ is $L$-periodic and with mean value zero with respect to $\tau$ and defined for $|z|\le K$ ($K$ is a positive constant depending only on $s',L$ and $\|I_0\|_{H^{s'}}$) and
$y\in\mathbb{R}^{n-1}$,
        \item  for $\Omega=\mathbb{R}/L\mathbb{Z}\times
\mathbb{R}^{n-1}_y$ 
$z\mapsto I(\tau,z,\mathbf{y}) \in C(]-K,K[; H^{s'}(\Omega)) \cap C^1(]-K,K[;H^{s'-2}(\Omega)).$
            \end{itemize}
%Then the remainder term $R$ in Eq.~(\ref{CCL})

   Let $\overline{U}_\epsilon$ be the
approximate solution of the isentropic Euler system deduced from a
solution of the KZK equation with the help of the correctors $v$, $w$, $v_1$, found by $I$ following the formulae of the derivation KZK-ansatz, ensuring the remainder term of the order of $\epsilon^\frac{5}{2}$. Then the function
$\overline U_\epsilon(x_1,\mathbf{x'},t)=\overline U_\epsilon(t-\frac{x_1}c,\epsilon
x_1,\sqrt {\epsilon}\mathbf{x'})$ given by  formula~(\ref{linU})
 is defined
in $$
  \mathbb{R}_t\times (\Omega_\epsilon= \{ |x_1|<\frac{K}{\epsilon}-ct\}\times
\mathbb{R}^{n-1}_{\mathbf{x'}})$$ and is smooth enough according to the above
procedure and the remainder term $\mathbf{R}$ in Eq.~(\ref{CL}) is in $[L_{\infty}((-K,K); L_2)]^2.$

Let us now consider  the solution of the Euler system~(\ref{CCL}) in
a cone (see  Fig.~\ref{Figura})
\begin{multline*} C(t)=\cup_{0<s<t}\{s\}\times
Q_\epsilon(s)=\\\{x=(x_1,\mathbf{x'}):\, |x_1|\le \frac{K}{\epsilon}-Ms, \; M\ge c , \;
\mathbf{x'} \in \mathbb{R}^{n-1}\}
\end{multline*} with
the initial data
\begin{eqnarray}
    &&(\bar{\rho}_\epsilon-\rho_\epsilon)|_{t=0} =0,
\quad (\mathbf{\bar{u}_\epsilon-u_\epsilon})|_{t=0} = 0.\label{bban}
\end{eqnarray}
\begin{figure}[th]
        \centering
        \psfrag{O}{$0$}
        \psfrag{x}{$\mathbf{x'}$}
        \psfrag{y}{$x_1$}
        \psfrag{t}{$t$}
        \psfrag{u}{$T$}
        \psfrag{r}{$\frac{K}{\epsilon}$}
        \psfrag{s}{$-\frac{K}{\epsilon}$}
        \psfrag{p}{slope $-M$}
        \psfrag{q}{slope $M$}
        \centerline{\includegraphics[width=7cm]{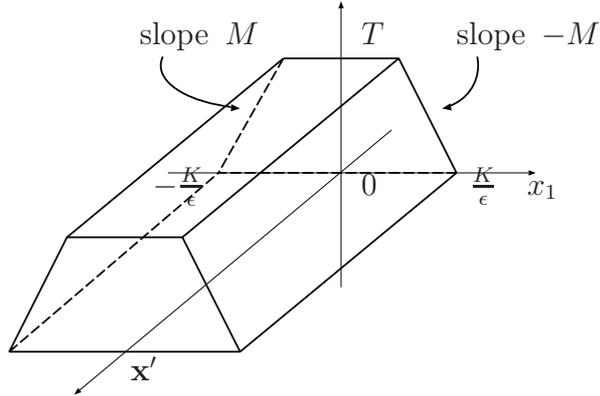}} 
        \caption{The cone  $C(T)$.}\label{Figura}
    \end{figure}
Consequently, there exists $T_0$  such
that for the time interval $0\leq t\leq \frac{T_0}{\epsilon}$ there
exists the classical solution $\mathbf{U_\epsilon}=(\rho_\epsilon,\mathbf{u_\epsilon} )$ of the
Euler system~(\ref{CCL}) in a cone
\begin{equation}\label{cone}
C(T)=\{0<t<T|T<\frac{T_0}{\epsilon}\}\times Q_\epsilon(t)
\end{equation}
 with $$\|\nabla. \mathbf{U_\epsilon}\|_{L_{\infty}([0,\frac{T_0}{\epsilon}[; H^{s'-5})}<\epsilon
    C \; \hbox{ for } \; s'>[\frac{n}{2}]+5.$$

Moreover, there exist positive constants $C_1$ and $C_2$  such that for any $\epsilon$ small
enough, the solutions
$\mathbf{\widetilde{U}_\epsilon}\stackrel{note}{=}(\rho_\epsilon,\rho_\epsilon \mathbf{u_\epsilon})$
and $\mathbf{\bar{
U}_\epsilon}\stackrel{note}{=}(\overline{\rho}_\epsilon,\overline{\rho}_\epsilon\mathbf{\overline{u}_\epsilon})$,
which were determined as above in  cone~(\ref{cone})
   with the same initial data~(\ref{bban}),
satisfy the estimate
\begin{equation}\label{elestim}
C_1\epsilon^\frac{7}{2}t\le\|\mathbf{\bar{U}_\epsilon - \widetilde{U}_\epsilon}\|^2_{L_2(Q_\epsilon(t))}\le \epsilon^5
e^{C_2\epsilon t}.
\end{equation}
\end{theorem}
Let now consider the viscous case.

For the viscous case we have
\begin{equation}\label{NCL}
\partial_t \mathbf{\widetilde{U}_{\epsilon}} +\nabla. \mathbf{F(\widetilde{U}_{\epsilon})}-\epsilon\nu \left[\begin{array}{c}0\\ \triangle \mathbf{u_{\epsilon}}\end{array}\right]=0
\end{equation}
for the exact system, and 
\begin{equation}\label{NCLA}
\partial_t \mathbf{\overline{U}_{\epsilon}} +\nabla. \mathbf{F(\overline{U}_{\epsilon})}-\epsilon\nu \left[\begin{array}{c}0\\ \triangle \mathbf{\overline{u}_{\epsilon}}\end{array}\right]=\epsilon^\frac{5}{3}\mathbf{R}
\end{equation}
for the approximate system.
\begin{theorem}\label{visqueux}
Suppose that  the initial data of the KZK Cauchy problem
$I_0(t,\mathbf{y})=I_0(t,\sqrt{\epsilon} \mathbf{x'})$ is such that \begin{enumerate}
    \item $I_0$ is $L$-periodic in $t$ and with  mean value zero, %\label{itT1}
    \item for fixed $t$, $I_0$ has the same sign for all $\mathbf{y}\in
    \mathbb{R}^{n-1}$, and for $t\in ]0,L[$ the sign changes, i.e. $I_0=0$, only for a finite number of times, %\label{itT2}
    \item $I_0(t,\mathbf{y})\in H^{s'}(\{t\ge0\}\times \mathbb{R}^{n-1})$ for
    $s'>\max\{6,[\frac{n}{2}]+1\}$,%\label{itT3}
    \item $I_0$ is sufficiently small %in the sense of Theorem~\ref{cauchy} 
    such that \begin{equation*}
 \|I_0\|_{H^{s'}}< \frac{\nu}{2c^2\rho_0}\frac{C_1(L)}{C_2(s')} \quad \hbox{(see \cite[p.20]{ARPAnal}),}
\end{equation*}
       and $I_0=\epsilon^\alpha \tilde{I}_0$, $\alpha\ge0$.%\label{itT4}
\end{enumerate}

Then there exists a
unique global  solution in time $\mathbf{\overline{U}_\epsilon}=(\bar{\rho}_\epsilon,
\mathbf{\bar{u}_\epsilon})$
 of the approximate system~(\ref{NCLA})
deduced from a solution of the KZK equation with the help
of correctors $v$, $\mathbf{w}$, $v_1$, found by $I$ following the formulae of the derivation KZK-ansatz, ensuring the remainder term of the order of $\epsilon^\frac{5}{2}$. The function $\mathbf{\overline{
U_\epsilon}}(x_1,\mathbf{x'},t)=\mathbf{\overline{U_\epsilon}}(x_1-ct,\epsilon x_1,\sqrt
{\epsilon}\mathbf{x'})$, given by  formula~(\ref{linU}), is defined in the
half space (see~\cite{ARPDer} for its regularity)
\begin{equation}\label{hs}
 \{x_1>0, \quad t>0, \quad \mathbf{x'}\in
\mathbb{R}^{n-1}\}.
\end{equation}
%  Moreover,
% we have
% \begin{equation}
% \bar{\rho}_\epsilon\in C\left([0,\infty[,
% H^{s'}\left(\Omega\right)\right) \cap
% C^1\left([0,\infty[;H^{s'-2}\left(\Omega\right)\right),\label{barro}
% \end{equation}
% \begin{equation}\bar{u}_\epsilon\in C\left([0,\infty[;
% H^{s'-2}\left(\Omega\right)\right) \cap
% C^1\left([0,\infty[;H^{s'-4}\left(\Omega\right)\right).\label{baru}
% \end{equation}

The Navier-Stokes system~(\ref{NCL}) in  the half space
 with
 initial data~(\ref{bban}) and following boundary conditions
$$    (\mathbf{\bar{u}_{\epsilon}-u_{\epsilon}})|_{x_1=0}=0,
$$%
with positive  first component of the velocity  
$u_{\epsilon,1}|_{x_1=0}
>0$ (i.e. at points where the fluid enters the domain) has the  additional boundary
condition
$$(\bar{\rho}_\epsilon-\rho_\epsilon)|_{x_1=0}=0.
$$
When $u_{\epsilon,1}|_{x_1=0} \le0$ there is no any boundary
condition for $\rho_\epsilon$.

Then there exists a constant $T_0>0$ such that for all $t<
    \frac{T_0}{\epsilon^{2+\alpha}}$ there exists a unique
    solution $\mathbf{U_\epsilon}=(\rho_\epsilon, \mathbf{u_\epsilon})$ of the Navier-Stokes system~(\ref{NCL}) with
    the same smoothness as $\mathbf{\overline
U_\epsilon}$. In addition, there exist positive constants $C_1>0$ and $C_2>0$ such that for all small enough
$\epsilon$
\begin{equation}\label{ocns}
C_1\epsilon^\frac{5}{2} \sqrt{t}\le\| \rho_\epsilon -\overline{\rho}_\epsilon \|_{L_2}+\|\rho_\epsilon
\mathbf{u_\epsilon}-\overline{\rho}_\epsilon \mathbf{\overline{u}_\epsilon}\|_{L_2} \le \epsilon^{\frac 52} e^{C_2\epsilon t}.
\end{equation}
Estimate~(\ref{ocns}) ensures that its left-hand side remains
smaller than the order of $\epsilon$ for any finite time
$$0<t<\frac{T}{\epsilon}\ln\frac{1}{\epsilon},$$ 
where $T$ is a positive constant and $T=O(1)$.
 
\end{theorem}

% Body of subsection.
% 
% Numbered formula should be written as
% \begin{equation}\label{key}                       % key is any name used
% a+b=2.                                         % to refer to the formula
% \end{equation}
% 
% The reference to such an equation should be as (\ref{key}).
% 
% A few formulae:
% \begin{gather}                       % key is any name used
% a+b=2,\label{key1}\\
% c+d=3.\label{key2}
% \end{gather}
% 
% 
% Aligned formulae:
% \begin{align}                       % key is any name used
% A={}&\int_0^1f(x)\mathrm{d}x,\notag\\
% &{}+a+b+c,\label{key3}\\
% B={}&3.\label{key4}
% \end{align}
% 
% A long formula:
% \begin{multline}
% D=a+b+c+e+\sum_{n=1}^\infty t_n\\
% +f+\int_0^1g(x)\mathrm{d}x+q+h+y^2\\
% -r-s-\sin(w).
% \end{multline}
%  
% % See http://mirror.ctan.org/macros/latex/required/amslatex/math/amsldoc.pdf
% % for more examples of AMS-LaTeX commands.
%  
% Citations are written as \cite{paper1}.           % paper1 is the name from
%                                                   % \bibitem commands below
% Fig. \ref{fig1} shows handling of figures.
% 
% \begin{figure}[t!]\centering
% %  \includegraphics[width=.43\textwidth]{Test_fig.epsilon}
%   \caption{This is test figure}\label{fig1}
% \end{figure}

% \section*{Acknowledgements}
% 
% Put acknowledgements in the last section, please do not use footnotes for that.


\begin {thebibliography}{99}

% \bibitem{paper1} Surname, I.I., Surname, I.I., 19??,
%             Title of the reference,
%             \emph{Journal}. Vol.\;{\bf ?}, pp.\;??--??.
% \bibitem{paper2} .....
\bibitem{Kuznetsov} Kuznetsov V.P., 1971, Equations of nonlinear acoustics, \emph{Sov. Phys. Acoust.}. Vol.\;{\bf 16},  pp.\;467--470.\bibitem{MCDKuper} McDonald B. E.,  Kuperman W.A., 1987, Time domain formulation for pulse propagation including nonlinear  behavior at a caustic, \emph{J. Acoust. Soc. Am.}. Vol.\;{\bf 81},  pp.\;1406--1417.
\bibitem{KhZabKuzn}  Bakhvalov N.S., 
Zhileikin Ya. M., Zabolotskaya E.A., 1987, Nonlinear Theory of Sound
Beams, \emph{American Institute of Physics}, New York.%  ({\em  Nelineinaya
%teoriya zvukovih puchkov}, Moscow ``Nauka", 1982).

\bibitem{ARPDer}  Rozanova-Pierrat A., 2009, On the derivation of the Khokhlov-Zabolotskaya-Kuznetsov (KZK)
  equation andvalidation of the KZK-approximation for viscous and non-viscous
  thermo-elastic media, \emph{Commun. Math. Sci.}. Vol.\;{\bf 7}, pp.\;679--718. 
    \bibitem{ARPAnal}  Rozanova-Pierrat A., 2008, Qualitative analysis of the Khokhlov–Zabolotskaya–Kuznetsov (KZK) equation,  \emph{Math. Mod. Meth. Appl. Sci.}. Vol.\;{\bf 18}, pp.\;781--812. 
\bibitem{Anna}  Rozanova A., 2007, Khokhlov-Zabolotskaya-Kuznetsov equation,
\emph{C. R. Acad. Sci. Paris, Ser. I}. Vol.\;{\bf 344}, pp.\;337--342. 
 \bibitem{Makarov} Makarov S., Ochmann M., 1997, Nonlinear and Thermoviscous Phenomena in Acoustics, Part II, \emph{Acustica}. Vol.\;{\bf 83}, pp.\;197--222.
 \bibitem{hamBlack}
  Hamilton M.F.,  Blackstock D.T., 1998,  Nonlinear Acoustics,
  \emph{Academic Press}, San Diego.
\bibitem{Jordan}  Jordan P.M., 2004, An analytical study of kuznetsov's equation: diffusive solitons,
 shock formation, and solution bifurcation,
\emph{Physics Letters A}. Vol.\;{\bf 326}, pp.\;77--84.




% \bibitem{kaltenbacher2011well}
% Kaltenbacher B., Lasiecka I., 2011,
%  Well-posedness of the Westervelt and the Kuznetsov equation with
%   nonhomogeneous neumann boundary conditions.
% {\em Discret. Contin. Dyn. S.}. Vol.\;{\bf  Issue Special}, pp.\;763--773.
% \bibitem{Dafermos}  Dafermos C., 2000, Hyperbolic Conservation Laws,
% \emph{Continuum Physics}, Springer-Verlag, 325.

\end{thebibliography}
\end{document}